\documentclass[a4paper,10pt,fleqn]{article}

\usepackage{theorem}
\usepackage{amssymb}
\usepackage{latexsym}

\theoremstyle{break}
\newtheorem{thm}{Theorem}
\newtheorem{prop}{Proposition}
\newtheorem{lem}{Lemma}
\newtheorem{de}{Definition}
\newtheorem{cor}{Corollary}
\newtheorem{nota}{Notation}
\newtheorem{conj}{Conjecture}

\title{\ A generalization of Chebyshev polynomials and non rooted posets}
\author{Masaya Tomie}
\date{tomie@math.tsukuba.ac.jp}

\begin{document}

\maketitle

In this paper we give a generalization of   Chebyshev polynomials 
and using this we describe the Mobius function of the 
generalized subword order 
derived from a poset 
$\{  a_{1}, \cdots a_{s} ,c \ | \ \ 
a_{i} < c \ for \ i = 1, \cdots s \}$, which contains an affirmative answer for the conjecture 
by Bj\"orner, Sagan and Vatter.(cf,\cite{bjorner5} \cite{sag})

\section{INTRODUCTION}

Bj\"orner was the first to determine  the M\"obius functions of 
factor orders and subword orders.
To determine the M\"obius functions, he used involutions, shellability, and 
generating functions. \cite{bjorner2}\cite{bjorner3}\cite{bjorner4}

Bj\"orner and Stanley found an interesting relation among 
the subword order derived from a two point set $\{ a,\ b \}$ , 
symmetric groups 
and  composition orders. \cite{bjorner6}

Factor orders,  subword orders, and 
generalized subword orders were studied 
in the context of M\"obius functions derived from word orders.

In \cite{sag} Sagan and Vatter gave a description of the M\"obius function of 
the generalized subword order derived from positive integers in two ways, 
namely 
the sign reversing involution and the discrete Morse theory.
More generally they gave a combinatorial description 
of the M\"obius functions derived from  rooted forests. 
And in \cite{bjorner5}\cite{sag}
 they gave a very  interesting conjecture which connects  
with the relation between a non-rooted forest $P_{2}$ as in 
Notation.{\rmfamily \ref{nonforest}} 
and Chebyshev polynomials.

\begin{conj}[\cite{bjorner5}\cite{sag}]\label{sagan-conj}

We put $P := \{ a,\ b,\ c, \ | \ a<c, \ b<c \ \}$, and consider the 
poset $P^{*}$ consisting of finite words of $P$ with its generalized subword 
order.
  Let $\mu$ be a M\"obius function of $P^{*}$. Suppose $0 \le i \le j$.

Then   $\mu(a^{i},c^{j})$  is the  
coefficient  of $ X^{j-i} $ in  $T_{i+j}(X)$.

Now we call $\{ T_{n}(X) \ | \ n \in \mathbb{N} \}$  Chebyshev polynomials of first kind.

\end{conj}

A series of Chebyshev polynomials  $\{ T(X) \ | \ n \in \mathbb{N} \}$ 
is a system of 
orthogonal polynomials and induces a special case of hypergeometric functions 
as a generalization of a binomial series. And this polynomial series 
is an example of the best approximation polynomials.

Not only in analysis, but in combinatorics, Chebyshev polynomials 
appear in permutation pattern avoidances \cite{chow} and 
Chebyshev posets, Chebyshev transformations defined by Hetyei which are 
related 
to $cd$-indeices, $f$-vectors and $h$-vectors respectively.

In this paper we give a  natural generalization of Chebyshev polynomials 
in the following way.

\begin{de}[generalized Chebyshev polynomials]

We define the polynomial $T^{s}_{k}(X) \ for \ s,k \in \mathbb{N}$ as follows:

(1) \ \ $T^{s}_{0}(X) = 1 , \ T^{s}_{1}(X) = (s-1)X$,

(2) \ \ $T^{s}_{k+2}(X) + T^{s}_{k}(X) = s X \cdot T^{s}_{k+1}(X)$.

Now the $T^{2}_{n}(X)$ are Chebyshev polynomials of first kind.
 And notice  $deg(T^{s}_{k}(X)) = k$.

\end{de}

Then, using generalized Chebyshev polynomials,  
 we generalize the conjecture as follows.

\begin{thm}\label{last1}
Let $P_{s}$ be a poset as Notation.{\rmfamily \ref{nonforest}} and 
$\mu$ be   the M\"obius function of $P_{s}^{*}$. Then 
for $0 \le m \le n$ ,   
$\mu (a_{1}^{m},c^{n})$ is the coefficient of $X^{n-m}$ in 
$T^{s}_{m+n}(x)$. 

\end{thm}

\section{PRELIMINERIES}

In this section, we give some basic definitions and notations used in this 
paper. 
For the basic definitions of posets and  M\"obius functions, see \cite{stan2} 
and for the 
 definitions of subword orders and generalized subword orders, 
see \cite{bjorner2} \cite{bjorner3} \cite{bjorner4} \cite{bjorner5} 
\cite{sag}.
First we recall a path in a poset $P$.\cite{stan2}

\begin{de}[\cite{stan2}]
Let $P$ be a poset. An arranged sequence $( \theta_{1}, \cdots , \theta_{r} )$ 
with  
$\theta_{i} \in P$ and $  \theta_{1}< \cdots < \theta_{r}, $  
is called a path of length $r-1$ .
We denote it  by $( \theta_{1} \rightarrow \cdots \rightarrow 
\theta_{r} )_{p}$

\end{de}

\begin{prop}[\cite{stan2}]\label{stanley}
Let $P$ be a locally finite poset, $\mu$ be a M\"obius function of $P$ and 
$C_{k}$ be the number of paths of $P$ of length $k$ respectively. Then we have 

$\mu (u,v) = \Sigma_{k} (-1)^{k} C_{k} $,
for all $u \le v \in P$.

\end{prop}

\begin{de}[\cite{sag}]
Let $P^{*}$ be the poset with the subword order derived from a poset 
$P$. We take 
$p_{1} , \cdots p_{k}$ and $q_{1} \cdots q_{l}$ from $P^{*}$. 
If $p_{1}  \cdots p_{k} \le q_{1} \cdots q_{l}$ as a subword order,  
we call $S(j_{1} , \cdots ,j_{k})$ 
($j_{1} < \cdots < j_{k}$), 
 an embedding of $p_{1}  \cdots p_{k}$ 
into $q_{1} \cdots q_{l}$ if $p_{i} \le q_{j_{i}} $ for  $1 \le i \le k$

And an embedding 
$S(j_{1} , \cdots ,j_{k})$ is called  the right most embedding of 
$p_{1}  \cdots p_{k}$ into $q_{1} \cdots q_{l}$, if for any embedding 
$S^{\prime}(j^{\prime}_{1} , \cdots , j^{\prime}_{k})$ , 
we have $j^{\prime}_{i} \le j_{i} $ for all $1 \le i \le k$.
\end{de}

\begin{nota}\label{nonforest}
In this paper we fix a poset $P_{s}$ for $s \in \mathbb{N}$ as follows.

$P_{s} := \{  a_{1}, \cdots a_{s} ,c \ | \ \ 
a_{i} < c , \ for \ i = 1, \cdots s \}         \}$

\end{nota}

\begin{de}
We define as follows.

Let $P_{s}^{*}$ be a poset with the generalized subword order derived from 
a poset $P_{s}$ 
as in Notation {\rmfamily \ref{nonforest}} and let $X$ be a set of the 
paths of $P$.

Put $Mob(X) := \Sigma_{k \ge 1} C_{k}$, where $C_{k}$ is the number of paths 
in $X$ whose length is $k$. Also we define 

$\{  a^{k}c^{l} \} := \{ p_{1} \cdots p_{k+l} \ | \ \ 
 \sharp \{ p_{1}, \cdots ,p_{k+l}  \} \cap \{ a_{1} , \cdots a_{s} \} = k , \ 
  \sharp \{ p_{1}, \cdots ,p_{k+l}  \} \cap \{ c \} = l \}$ for $k,l \in 
  \mathbb{N} \cup \{ 0 \}$, 

$< p_{1} \cdots p_{k} , \{ a^{l}c^{m} \} > := 
\{ q_{1} \cdots q_{l+m} \in  \{ a^{l}c^{m} \} \ | \ \ 
 p_{1} \cdots p_{k} \le   q_{1} \cdots q_{l+m}  \}
$ for $k,l \in  \mathbb{N} \cup \{ 0 \}$, and 

$Pat \{ p_{1} \cdots p_{k} , q_{1} \cdots q_{l}    \}
:= \{  ( p_{1} \cdots p_{k} \rightarrow \theta_{1} \rightarrow , \cdots 
\rightarrow \theta_{r} \rightarrow q_{1} \cdots q_{l} )_{p}   
\ | \ p_{1} \cdots p_{k} < \theta_{1}< \cdots < \theta_{r}
 < q_{1} \cdots q_{l} , \ |\theta_{i}| 
=l \}$ respectively.
Here $|\theta|$ is the number of letters of $\theta$.

\end{de}

\begin{de}[generalized Chebyshev polynomial]

We define the polynomial $T^{s}_{k}(X) \ s,k \in \mathbb{N}$ as follows:

(1) \ \ $T^{s}_{0}(X) = 1 , \ T^{s}_{1}(X) = (s-1)X$,

(2) \ \ $T^{s}_{k+2}(X) + T^{s}_{k}(X) = s X \cdot T^{s}_{k+1}(X)$.

Now the $\{ T^{2}_{n}(X) \ | \ n \in \mathbb{N} \}$  are Chebyshev 
polynomials of first kind.
 And notice  $deg(T^{s}_{k}(X)) = k$.

\end{de}

Here we give a simple expression of the generalized Chebyshev polynomials.

\begin{prop}

For $s, n \in \mathbb{N}$, we have

\[
T_{n}^{s}(X) = \Sigma_{m \le n, \ n-m : even}
(-1)^{(n-m)/2} \Bigl(
(\begin{array}{c}
(n+m)/2 \\
(n-m)/2
\end{array}
)
\cdot s^{m}
-
(\begin{array}{c}
(n+m)/2 -1 \\
(n-m)/2
\end{array}
)
\cdot s^{m-1}
\Bigr)
\cdot
X^{m}
\]

\end{prop}

 \textbf{PROOF}

We show the above formula 
by induction. It is easy to see $T_{0}^{s}(X) = 1, \ T_{1}^{s}(X)
 = (s-1)X$.
Now we have 

$T_{n}^{s} + T_{n+2}^{s} $

\[
 = \Sigma_{m \le n, \ n-m : even}
(-1)^{(n-m)/2} \Bigl(
(\begin{array}{c}
(n+m)/2 \\
(n-m)/2
\end{array}
)
\cdot s^{m}
-
(\begin{array}{c}
(n+m)/2 -1 \\
(n-m)/2
\end{array}
)
\cdot s^{m-1}
\Bigr)
\cdot
X^{m}
\]

\[
+ \Sigma_{m \le n+2, \ n+2-m : even}
(-1)^{(n+2-m)/2} \Bigl(
(\begin{array}{c}
(n+2+m)/2 \\
(n+2-m)/2
\end{array}
)
\cdot s^{m}
-
(\begin{array}{c}
(n+2+m)/2 -1 \\
(n+2-m)/2
\end{array}
)
\cdot s^{m-1}
\Bigr)
\cdot
X^{m}
\]

\[
= 
\Bigl(
(\begin{array}{c}
n+2 \\
0
\end{array})
\cdot
s^{n+2}
-
(\begin{array}{c}
n+2 \\
0
\end{array})
\cdot
s^{n+1}
\Bigr) \cdot X^{n+2}
\]

\[
+
\Sigma_{m \le n, \ n-m : even}
(-1)^{(n-m)/2 +1} \Bigl(
(\begin{array}{c}
(n+m)/2 \\
(n-m)/2 +1
\end{array}
)
\cdot s^{m}
-
(\begin{array}{c}
(n+m)/2 -1 \\
(n-m)/2 +1
\end{array}
)
\cdot s^{m-1}
\Bigr)
\cdot
X^{m}
\]

\[
= s^{n+1}(s-1)X^{n+2}
\]

\[
+
\Sigma_{m \le n-1, \ n-1-m : even}
(-1)^{(n-m-1)/2 +1} \Bigl(
(\begin{array}{c}
(n+m+1)/2 \\
(n-m-1)/2 +1
\end{array}
)
\cdot s^{m+1}
-
(\begin{array}{c}
(n+m+1)/2 -1 \\
(n-m-1)/2 +1
\end{array}
)
\cdot s^{m}
\Bigr)
\cdot
X^{m+1}
\]

\[
=sX \biggl(
s^{n}(s-1)X^{n+1}
\]

\[
+
\Sigma_{m \le n-1, \ n-1-m : even}
(-1)^{(n+1-m)/2 } \Bigl(
(\begin{array}{c}
(n+1+m)/2 \\
(n+1-m)/2 
\end{array}
)
\cdot s^{m}
-
(\begin{array}{c}
(n+1+m)/2 -1 \\
(n+1-m)/2 
\end{array}
)
\cdot s^{m-1}
\Bigr)
\cdot
X^{m}
\Biggr)
\]

$= sX \cdot T_{n+1}^{s}$.

Hence we obtain the derived result. \ \ \ \ \ \ \ \ \ \ \  \ \ \ \ \ \
$\Box$

\section{MAIN RESULTS}

In this section, we give a proof of Theorem {\rmfamily \ref{last1}}

\begin{lem}\label{mobius}

Let $P$ be a finite poset and we take an element $x \in P$.  We put as follows:

$P_{\le x} := \{ y \ | \ y \le x  \}$ , 
$\widehat{P_{\le x}}
:= \{ (\theta_{1} \rightarrow  \cdots 
\rightarrow \theta_{r-1} \rightarrow x)_{p} \ | \ \theta_{i} \in P  \}$ , 
$P_{\ge x} := \{ y \ | \ y \ge x  \}$ , 
$\widehat{P_{\ge x}} := \{ 
x \rightarrow \theta_{1} \rightarrow \cdots \rightarrow \theta_{r-1} )_{p}
 \ | \ \theta_{i} \in P\}$
and 
$\mathfrak{P}_{x} := \{ ( \cdots \rightarrow 
\tau_{r} \rightarrow x \rightarrow \sigma_{1} \rightarrow \cdots )  \ | \ 
\tau_{i} \le x , \ \sigma_{i} \ge x \   \}$.
Now a path 
$(x) \in \widehat{P_{\le x}}, \ \widehat{P_{\ge x}}, \ \mathfrak{P}_{x}$.
Then we have $Mob \mathfrak{P}_{x}  =  Mob \widehat{P_{\le x}} \ 
 Mob \widehat{P_{\ge x}} $.
 
 \end{lem}

 \textbf{PROOF}

Notice that a path which passes through $x$ splits into the two paths, one 
starts from 
$x$ and the other one ends $x$. From that we obtain the derived result.
\ \ \ \ \ $\Box$

\begin{lem}\label{firstcombi}

For $m,n,p,q \in \mathbb{N} \cup \{ 0 \}$ such that 
$0 \le m \le n  , \ 0 \le p \le m , \ 0 \le q \le n 
$, 
 we take $p_{1} \cdots p_{m} , \tilde{p_{1}} \cdots \tilde{p_{m}}
 \in \{ a^{m-p}c^{p} \} $.
Then we have 

$\sharp <p_{1} \cdots p_{m} , \{ a^{n-q}c^{q} \} > 
= \sharp < \tilde{p_{1}} \cdots \tilde{p_{m}} , \{ a^{n-q}c^{q} \} >$.

\end{lem}

\textbf{PROOF}

\underline{Claim 1} \ \ 
We have $\sharp <p_{1} \cdots 
\displaystyle\overbrace{a_{x}}^{i-th}  \cdots  p_{m} , \{ a^{n-q}c^{q} \} > 
= \sharp <p_{1} \cdots 
\displaystyle\overbrace{{a_{y}}}^{i-th} \cdots 
p_{m} , \{ a^{n-q}c^{q} \} >$. 

(Proof of claim1)

We take $\forall q_{1} \cdots q_{n} \in 
<p_{1} \cdots \displaystyle\overbrace{a_{x}}^{i-th} 
  \cdots  p_{m} , \{ a^{n-q}c^{q} \} >$. And we  consider 
the right most embedding  into 
$q_{1} \cdots q_{n}$. Notice that the right most embedding is unique.
Here we put $S(j_{1},j_{2}, \cdots , j_{m} )$ as the 
 right most embedding 
$p_{1} \cdots \displaystyle\overbrace{a_{x}}^{i-th}
 \cdots  p_{m}$ into $q_{1} \cdots q_{n}$ 
.

Now we define the map $\Phi$ as follows.

$\Phi$ \ \ \ \ \ \ \ \ \ \ \ \ \ \ 
$<p_{1} \cdots \displaystyle\overbrace{a_{x}}^{i-th}
 \cdots  p_{m} , \{ a^{n-q}c^{q} \} >
 \longrightarrow 
 <p_{1} \cdots \displaystyle\overbrace{a_{y}}^{i-th}
  \cdots p_{m} , \{ a^{n-q}c^{q} \} >$

$\Phi (q_{1} q_{2}, \cdots 
\displaystyle\overbrace{a_{x}}^{j_{i}-th} 
a_{x_{1}} \cdots a_{x_{t}}   
\displaystyle\overbrace{p_{i+1}}^{j_{i+1}-th} 
  \cdots q_{n} ) = 
q_{1} q_{2} \cdots 
\displaystyle\overbrace{a_{y}}^{j_{i}-th} 
a_{y_{1}} \cdots a_{y_{t}} 
\displaystyle\overbrace{p_{i+1}}^{j_{i+1}-th}
 \cdots  q_{n} $,

Here we put $a_{y_{k}} = a_{x_{k}+y-x}$, and $a_{k+s} = a_{k}$.

It is easy to see  the right most embedding of 
$p_{1} \cdots \displaystyle\overbrace{a_{y}}^{i-th} \cdots p_{m} $
 into $\Phi (q_{1} \cdots q_{n})$ is $S(j_{1},j_{2}, \cdots , j_{m} )$.
And by the construction of $\Phi$, we can easily define the inverse map of 
$\Phi$.
Hence we prove this claim.
 
 \underline{Claim 2} \ \ 
We have $\sharp <p_{1} \cdots 
\displaystyle\overbrace{a_{x}}^{i-th}  
\displaystyle\overbrace{c}^{(i+1)-th}
 \cdots  p_{m} , \{ a^{n-q}c^{q} \} > 
= \sharp <p_{1} \cdots 
\displaystyle\overbrace{c}^{i-th} 
\displaystyle\overbrace{a_{x}}^{(i+1)-th} \cdots 
p_{m} , \{ a^{n-q}c^{q} \} >$. 

(Proof of claim2)

We take $\forall q_{1} \cdots q_{n} \in 
<p_{1} \cdots \displaystyle\overbrace{a_{x}}^{i-th} 
\displaystyle\overbrace{c}^{(i+1)-th} 
  \cdots  p_{m} , \{ a^{n-q}c^{q} \} >$ and put 
  $S(j_{1},j_{2}, \cdots ,  j_{m} )$ as  
the right most embedding 
$p_{1} \cdots \displaystyle\overbrace{a_{x}}^{i-th} 
\displaystyle\overbrace{c}^{(i+1)-th} 
  \cdots  p_{m}$ into 
$q_{1} \cdots q_{n}$. 

Now we define the map $\Phi$ as follows.

$\Phi$ \ \ \ \ \ \ \ \ \ \ \ \ \ \ 
$<p_{1} \cdots \displaystyle\overbrace{a_{x}}^{i-th} 
\displaystyle\overbrace{c}^{(i+1)-th}
 \cdots  p_{m} , \{ a^{n-q}c^{q} \} >
 \longrightarrow 
 <p_{1} \cdots \displaystyle\overbrace{c}^{i-th}
  \displaystyle\overbrace{a_{x}}^{(i+1)-th}
  \cdots p_{m} , \{ a^{n-q}c^{q} \} >$,
  
  $\Phi (q_{1} \cdots \displaystyle\underbrace{
  \displaystyle\overbrace{q_{j_{i}}}^{=a_{x}}
  \cdots }_{A} 
  \displaystyle\underbrace{
  \displaystyle\overbrace{q_{j_{i+1}}}^{=c} \cdots }_{B}
   q_{j_{i+2}} \cdots q_{n}) = 
  q_{1} \cdots \displaystyle\underbrace{
  \displaystyle\overbrace{q_{j_{i+1}}}^{=c}
  \cdots }_{B} 
  \displaystyle\underbrace{
  \displaystyle\overbrace{q_{j_{i}}}^{=a_{x}} \cdots }_{A}
   q_{j_{i+2}} \cdots q_{n}$.
  
  Here the right most embedding of $p_{1} \cdots 
\displaystyle\overbrace{c}^{i-th} 
\displaystyle\overbrace{a_{x}}^{(i+1)-th} 
\cdots p_{m} $ into 
$q_{1} \cdots \displaystyle\underbrace{
  \displaystyle\overbrace{q_{j_{i+1}}}^{=c}
  \cdots }_{B} 
  \displaystyle\underbrace{
  \displaystyle\overbrace{q_{j_{i}}}^{=a} \cdots }_{A}
   q_{j_{i+2}} \cdots q_{n}$ is  
$S(j_{1} \cdots j_{i}, j_{i+2}+j_{i}-j_{i+1} , j_{i+2} \cdots j_{m} )$.
 By the construction, all of the elements of  
 
 $<p_{1} \cdots 
\displaystyle\overbrace{a_{x}}^{i-th}  
\displaystyle\overbrace{c}^{(i+1)-th}
 \cdots  p_{m} , \{ a^{n-q}c^{q} \} >$ 
 whose right most embedding 
 are  $S(j_{1},j_{2}, \cdots , j_{m} )$, 
 have one to one correspondence to the elements of  
 $<p_{1} \cdots 
\displaystyle\overbrace{c}^{i-th} 
\displaystyle\overbrace{a_{x}}^{(i+1)-th} \cdots 
p_{m} , \{ a^{n-q}c^{q} \} >$ 
whose 
 right most embedding are  
 $S(j_{1} \cdots j_{i}, j_{i+2}+j_{i}-j_{i+1} , j_{i+2} \cdots j_{m} )$.
 Hence the  $\Phi$ is bijeciton.
Therefore we have this claim2.
By these claims we obtain the derived result.  \ \ \ \ \ \ \ $\Box$ 

From Lemma {\rmfamily \ref{firstcombi}} ,
 if $p_{1} \cdots p_{m} \in \{ a^{m-p}c^{p}   \}$, then 

$\sharp <p_{1} \cdots p_{m} , \{ a^{n-q}c^{q} \} > 
= \sharp <
\displaystyle\underbrace{a_{1} \cdots a_{1}}_{(m-p) times}
\displaystyle\underbrace{c \cdots c}_{p-times} , 
\{ a^{n-q}c^{q} \}>$. So we denote the number as 
$M((m,p),(n,q))$ for all 
$0 \le m \le n , \ 0 \le p \le m , \ 0 \le q \le n$.

\begin{lem}\label{boole}

Let $k, l \in  \mathbb{N} \cup \{ 0 \}  , 0 \le k \le l$ 
and  $p_{1} \cdots p_{l} \in 
\{ a^{l-k}c^{k} \}$, then we have 
 $[p_{1} \cdots p_{l} , c^{l}] \simeq 
B_{l-k}$. Now $B_{l-k}$ is a Boolean algebra of rank $l-k$.

\end{lem}

\begin{lem}\label{mobpatk}

For $m,n,p,k \in \mathbb{N} \cup \{0 \}$ such that 
$0 \le m \le n , \ 0 \le p \le m$, 
we take $p_{1} \cdots p_{m} \in \{ a^{m-p}c^{p} \}$, 
then the number of paths in $Pat \{ p_{1} \cdots p_{m} , c^{n} \}$ 
whose length are $k$ equals 
to the number of paths in  $Pat \{
\displaystyle\underbrace{a_{1} \cdots a_{1}}_{(m-p) times}
\displaystyle\underbrace{c \cdots c}_{p-times}  , c^{n} \}$.

\end{lem}

\textbf{PROOF}

Notice that if we take $q_{1} \cdots q_{n} \in \{a^{n-q} c^{q} \}$, 
then the number of length $l$ 
paths from $q_{1} \cdots q_{n}$  to $c^{n}$ 
equals to the number of length $l$ paths from 
$\displaystyle\underbrace{a_{1} \cdots a_{1}}_{(n-q) times}
\displaystyle\underbrace{c \cdots c}_{q-times} $ to $c^{n}$.

Hence we have

$ \sharp \{ p_{1} \cdots p_{m} \rightarrow  \theta_{1} 
\rightarrow \cdots \rightarrow \theta_{k} =c^{n} \ | \ |\theta_{i}| = n  
\ for \ i = 1, \cdots k 
 \}$

$= \Sigma_{p \le r \le n} M((m,p),(n,r)) \sharp \{  
\displaystyle\underbrace{a_{1} \cdots a_{1}}_{(n-r) times}
\displaystyle\underbrace{c \cdots c}_{r-times}  \rightarrow 
\tau_{1} \rightarrow \cdots \rightarrow \tau_{k-1} = c^{n}  \} $.
 \ \ \ \ \ \ \ \ (By Lemma {\rmfamily \ref{boole}})

Hence we obtain the derived result.

\ \ \ \ \ \ \ \ \ \ $\Box$

\begin{lem}

For $m,n,p \in \mathbb{N} \cup \{0 \}$, such that 
$0 \le m \le n , \ 0 \le p \le m$, 
we take $p_{1}, \cdots p_{m} \in \{ a^{m-p}c^{p} \}$. Then we have 

\[
MobPat \{  p_{1} \cdots p_{m} , c^{n} \} 
= 
\left\{
\begin{array}{rl}
- \Sigma_{i=0}^{n-p} (-1)^{n-p-i} M((m,p,),(n,p+i)) & if m < n \\
(-1)^{n-m} & if m=n
\end{array}
\right.
\]

\end{lem}

\textbf{PROOF}

From Lemma {\rmfamily \ref{mobpatk}} we have 
$MobPat \{  p_{1} \cdots p_{m} , c^{n} \} = 
MobPat \{ \displaystyle\underbrace{a_{1} \cdots a_{1}}_{(m-p) times}
\displaystyle\underbrace{c \cdots c}_{p-times} , c^{n} \}$.

Then we have 

$MobPat \{ \displaystyle\underbrace{a_{1} \cdots a_{1}}_{(m-p) times}
\displaystyle\underbrace{c \cdots c}_{p-times} , c^{n} \}$

$= (-1) \Sigma_{i=0}^{n-p} M((m,p),(n,p+i)) \mu 
(\displaystyle\underbrace{a_{1} \cdots a_{1}}_{(n-p-i) times}
\displaystyle\underbrace{c \cdots c}_{(p+i) times} , c^{n})$
\ \ \ \ \ \ \ \ \ \ (By Proposition {\rmfamily \ref{stanley}})

$= (-1) \Sigma_{i=0}^{n-p} M((m,p),(n,p+i)) (-1)^{n-p-i}$.

Hence we obtain the derived result. \ \ \ \ \ \ \ \ \ \ $\Box$

\begin{lem}\label{M-first}

For $m,n,p ,q \in \mathbb{N} \cup \{0 \}$ such that 
$1 \le m \le n , \ 1 \le p \le m, \ 1 \le q \le n, \ 1 \le p \le q$, 
 we have 

$M((m,p),(n,q)) = \Sigma_{i=0}^{n-m} M((m-1,p-1),(n-1-i,q-1)) \cdot s^{i}$.

\end{lem}

\textbf{PROOF}

We have 

the left hand side 

$= \sharp \{ x_{1} \cdots x_{n} \in  \{ a^{n-q}c^{q} \} \ | \   a_{1}^{m-p}c^{p}  
\le x_{1} \cdots x_{n} \}$ \ \ \ \ \ \ \ \ (By Lemma  
{\rmfamily \ref{firstcombi}})

$= \sharp \{  x_{1} \cdots x_{n} \in  \{ a^{n-q}c^{q} \} \ | \  
a_{1}^{m-p}c^{p}  \le x_{1} \cdots x_{n} ,  x_{n}=c       \}$

 $+ \sharp \{  x_{1} \cdots x_{n} \in  \{ a^{n-q}c^{q} \} \ | \  
 a_{1}^{m-p}c^{p} 
  \le x_{1} \cdots x_{n} ,  x_{n-1}=c , x_{n} \neq c     \}$
  
  $\cdots$
  
$+ \sharp \{  x_{1} \cdots x_{n} \in  \{ a^{n-q}c^{q} \} \ | \  
a_{1}^{m-p}c^{p} 
  \le x_{1} \cdots x_{n} ,  x_{n-i}=c , x_{n-i+k} \neq c  (1 \le k \le i)   \}$  
$\vdots$

$+ \sharp \{  x_{1} \cdots x_{n} \in  \{ a^{n-q}c^{q} \} \ | \ 
 a_{1}^{m-p}c^{p} 
  \le x_{1} \cdots x_{n} ,  x_{m}=c , x_{m+k} \neq c  (1 \le k \le n-m)   \}$

$= \sharp \{ \displaystyle\underbrace{\cdots}_{A}   \ | 
\ a_{1}^{m-p}c^{p-1} \le A \ A \in \{ a^{n-q}c^{q-1} \}  \}$

$+ \sharp \{ \displaystyle\underbrace{\cdots}_{A}  
\displaystyle\underbrace{c}_{n-1}  
 \ | 
\ a_{1}^{m-p}c^{p-1} \le A \ A \in \{  a^{n-q-1}c^{q-1} \}  \} \cdot s^{1}$

$\vdots$

$+ \sharp \{ \displaystyle\underbrace{\cdots}_{A}  
\displaystyle\underbrace{c}_{n-i}  
 \ | 
\ a_{1}^{m-p}c^{p-1} \le A \ A \in \{ a^{n-q-i}c^{q-1} \} \} \cdot s^{i}$

$\vdots$

$+ \sharp \{ \displaystyle\underbrace{\cdots}_{A}  
\displaystyle\underbrace{c}_{n-i}  
 \ | 
\ a_{1}^{m-p}c^{p-1} \le A \ A \in \{  a^{m-q}c^{q-1} \} \} \cdot s^{n-m}$

(In case of $n-q-i <0$ , we recognize 
$\sharp \{ \displaystyle\underbrace{\cdots}_{A}  
\displaystyle\underbrace{c}_{n-i}  
 \ | 
\ a^{m-p}c^{p-1} \le A \ A \in a^{n-q-i}c^{q-1} \}$ as $0$.)

$=M((m-1,p-1),(n-1,q-1)) \cdot s^{0} + M((m-1,p-1),(n-2,q-1)) \cdot s^{1}
+ \cdots $

$+M((m-1,p-1),(n-1-i,q-1)) \cdot s^{i}  \cdots 
M((m-1,p-1),(m-1,q-1)) \cdot s^{n-m}$

Hence we obtain the derived result. \ \ \ \ \ \ \ \ \ \ $\Box$

\begin{lem}

For $m,n,i \in \mathbb{N} \cup \{0 \}$ such that 
$i \le m \le n , \ i \le p \le q, \  i \le p \le m , \ i \le q \le n$,
we have

\[ M((m,p),(n,q)) = \Sigma_{k=0}^{n-m} M((m-i,p-i),(n-i-k,q-i)) \cdot 
s^{k} \cdot 
( \begin{array}{c}
i+k-1  \\ 
i-1 
\end{array} ).
\]

\end{lem}

\textbf{PROOF}

In  case of $i=1$, it is shown by Lemma {\rmfamily \ref{M-first}}.

We show the above formula by induction. 
We suppose this lemma holds for $i-1$.
Now we see

\[
M((m,p),(n,q)) = 
\Sigma_{k=0}^{n-m} M((m-i+1,p-i+1),(n-i-k+1,q-i+1)) \cdot s^{k}
 (\begin{array}{c}
 i+k-2 \\
 i-2
 \end{array})
 \]
 
 \[
 =\Sigma_{k=0}^{n-m} 
 (\Sigma_{l=0}^{n-m}M((m-i,p-i),(n-i-k-l,q-i)) \cdot 2^{l}) \cdot s^{k}
(\begin{array}{c}
 i+k-2 \\
 i-2
 \end{array})
 \]

\[
= \Sigma_{k,l =0}^{n-m} M((m-i,p-i),(n-i-(k+l),q-i)) \cdot s^{k+l}
(\begin{array}{c}
 i+k-2 \\
 i-2
 \end{array})
 \]

\[
= \Sigma_{k+l =0}^{n-i} M((m-i,p-i),(n-i-(k+l),q-i)) \cdot s^{k+l}
(\begin{array}{c}
 i+k-2 \\
 i-2
 \end{array})
 \]

\[
= \Sigma_{\alpha =0}^{n-i} M((m-i,p-i),(n-i-\alpha,q-i))
( \Sigma_{j=0}^{\alpha} 
(\begin{array}{c}
 i+j-2 \\
 i-2
 \end{array})
) \cdot s^{\alpha} .
\]

Now we remark the following formula.

\[
\Sigma_{x=0}^{\alpha} 
(\begin{array}{c}
 x+i \\
 i
 \end{array})
 = 
(\begin{array}{c}
 i+\alpha+1 \\
 i+1
 \end{array}),
\]

hence we have 

\[
= \Sigma_{\alpha =0}^{n-i} M((m-i,p-i),(n-i-\alpha,q-i))
( \Sigma_{j=0}^{\alpha} 
(\begin{array}{c}
 i+j-2 \\
 i-2
 \end{array})
) \cdot s^{\alpha} 
\]

\[
= \Sigma_{\alpha =0}^{n-i} M((m-i,p-i),(n-i-\alpha,q-i))\cdot s^{\alpha} 
(\begin{array}{c}
 n+\alpha-1 \\
 \alpha-1
 \end{array}).
\]

\ \ \ \ \ \ \ \ \ \ $\Box$

\begin{lem}

For $m,n,p,q \in \mathbb{N} $,  
$1 \le m \le n , \ 1 \le p \le q, 1 \le  p \le m , \ 1 \le q \le n,$
we have

\[ M((m,p),(n,q)) = \Sigma_{k=0}^{n-m} M((m-p,0),(n-p-k,q-p)) \cdot 
s^{k} \cdot 
( \begin{array}{c}
p+k-1  \\ 
p-1 
\end{array} ).
\]

\end{lem}

\begin{lem}

For $ 1 \le \alpha \le \beta$, 
we have 

$\Sigma_{i=0}^{\beta} M((\alpha,0),(\beta,i)) \cdot (-1)^{i} = 0$.

\end{lem}

\textbf{PROOF}

We give a combinatorial proof.
We put $\widehat{M}_{i} :=
< \displaystyle\underbrace{a_{1} \cdots a_{1}}_{\alpha-times}, 
\{ a^{\beta-i}c^{i} \} >$, 
$\widehat{M} := \biguplus_{0 \le i \le \beta} \widehat{M}_{i}$,
 $\widehat{M}_{ev} := \biguplus_{0 \le i \le \beta \ i; \ even} 
 \widehat{M}_{i}$ and 
 $\widehat{M}_{odd} := \biguplus_{0 \le i \le \beta \ i; \ odd} 
 \widehat{M}_{i}$. 

Then we have 
$\Sigma_{i=0}^{\beta} M((\alpha,0),(\beta,i)) \cdot (-1)^{i} = 
\sharp \widehat{M}_{ev}  -   \sharp \widehat{M}_{odd}$.

We consider the map $\Psi$ as follows.

$\Psi$ \ \ \ \ \ \ \ \ \ \ \ \ \ \ $\widehat{M} \longrightarrow \widehat{M}$

$\Psi(\displaystyle\underbrace{\cdots}_{A}  a_{1} a_{x_{1}} 
\cdots a_{x_{t}}) = 
\displaystyle\underbrace{\cdots}_{A}  c a_{x_{1}} 
\cdots a_{x_{t}} )$

$\Psi(\displaystyle\underbrace{\cdots}_{A}  c 
a_{x_{1}}  
\cdots a_{x_{t}}) = 
\displaystyle\underbrace{\cdots}_{A}  a_{1} 
a_{x_{1}}  
\cdots a_{x_{t}}$

$\Psi(\displaystyle\underbrace{\cdots}_{A}  a_{1}) = 
\displaystyle\underbrace{\cdots}_{A}  c$

$\Psi(\displaystyle\underbrace{\cdots}_{A}  c) = 
\displaystyle\underbrace{\cdots}_{A}  a_{1}$

Here $\Psi$ changes $a_{1}$ into $c$ and $c$ into $a_{1}$ 
which appears right most position 
of each elements. Since $\alpha$ not being $0$, 
each element of $\widehat{M}$ contains $a_{1}$ or $c$. From that the map 
$\Psi$ is well-defined. Therefore obviously $\Psi^{-1} = \Psi$ and 
$\Psi (\widehat{M}_{ev}) = (\widehat{M}_{odd}) \ 
\Psi (\widehat{M}_{odd}) = \widehat{M}_{ev}$. Hence $\Psi$ is a bijection and 
$\sharp \widehat{M}_{ev} = \sharp \widehat{M}_{odd}$.
Hence we obtain the derived result.
\ \ \ \ \ \ \ \ \ \ $\Box$

\begin{lem}\label{import}

For $m,n,p \in \mathbb{N} \cup \{0 \}$,  
$0 \le m < n , \ 0 \le p < m$, 
 we have 

$p_{1} \cdots p_{m} \in \{ a^{m-p}c^{p}   \}  \Longrightarrow MobPat 
\{ p_{1}  \cdots p_{m}  , c^{n}  \} =0$.

\end{lem}

\textbf{PROOF}

We have 

$MobPat \{ p_{1} \cdots p_{m}  , c^{n}    \}
= - \Sigma_{i=0}^{n-p} M((m,p),(n,p+i)) \cdot (-1)^{n-p-i}$

\[
= -\Sigma_{i=0}^{n-p} (-1)^{n-p-i} \Sigma_{k=0}^{n-m}
M((m-p,0),(n-p-k,i)) \cdot s^{k}
(\begin{array}{c}
p+k-1 \\
p-1
\end{array})
\]

\[
= (-1)^{n-p-1} \Sigma_{i=0}^{n-p} \Sigma_{k=0}^{n-m} M((m-p,0),(n-p-k,i))
(-1)^{i} \cdot s^{k} 
(\begin{array}{c}
p+k-1 \\
p-1
\end{array})
\]

\[
= (-1)^{n-p-1} \Sigma_{k=0}^{n-m}
 \{ \displaystyle\underbrace
 {\Sigma_{i=0}^{n-p-k}M((m-p,0),(n-p-k,i))(-1)^{i}}_{=0}  \}
 \cdot s^{k}
(\begin{array}{c}
p+k-1 \\
p-1
\end{array})
\]

$=0$.   Hence we obtain the derived result. \ \ \ \ \ \  \ \ \ \ $\Box$

\begin{lem}

For $m,n \in \mathbb{N}, \ 1 \le m \le n$,
 we put

$P:= \{ a_{1}^{m} \rightarrow \tau_{1} \rightarrow 
\cdots \rightarrow c^{m} \rightarrow 
\theta^{k_{1}}_{1} \cdots \rightarrow 
c^{k_{1}} \rightarrow 
\theta^{k_{2}}_{1} \cdots \rightarrow 
 c^{k_{2}} \rightarrow \cdots \rightarrow c^{k_{r}}  \ | \ 
m < k_{1} < \cdots < k_{r}=n, \ | \tau_{i} | = m, \ 
|\theta^{k_{i}}_{j}| = k_{i}    \  \}$, 
and $p_{1} \cdots p_{m} \in \{ a^{m}  \ 
\}$.

Then we have 
$\mu (p_{1} \cdots p_{m}, c^{n}) = 
\mu (a_{1}^{m},c^{n}) = Mob(P) = (-1)^{m} \mu (c^{m},c^{n})$.

\end{lem}

\textbf{PROOF}

We have 

$\mu (p_{1} \cdots p_{m}, c^{n})$
$= Mob (\{ (p_{1} \cdots p_{m} \rightarrow \theta_{1} \rightarrow 
\cdots \theta_{r} \rightarrow c^{n}) \ | \  
p_{1} \cdots p_{m} < \theta_{1} < \cdots \theta_{r} < c^{n}
     \})$.

Now we put 

$X^{l_{1}}_{l_{2}}$

$ := \{  (
p_{1} \cdots p_{m} \rightarrow \cdots \theta_{r} \rightarrow 
\tau_{1} \rightarrow \cdots \tau_{s} \rightarrow c^{l_{2}} \rightarrow 
\sigma^{k_{1}}_{1} \rightarrow \cdots \rightarrow c^{k_{1}} \rightarrow 
 \sigma^{k_{2}}_{1} \rightarrow \cdots \rightarrow c^{k_{2}} \cdots 
 \rightarrow c^{n} ) \ | \ 
 | \theta_{r} | = l_{1}, \ \theta_{r} \neq c^{l_{1}}, \ 
 |\tau_{1}|, \cdots |\tau_{s}| = l_{2} , \ |\sigma^{k_{i}}_{j}| = k_{i}    \}$.

Then we have 

$\{ (p_{1} \cdots p_{m} \rightarrow \theta_{1} \rightarrow 
\cdots \theta_{r} \rightarrow c^{n}) \ | \  
p_{1} \cdots p_{m} < \theta_{1} < \cdots \theta_{r} < c^{n}
     \}
 $

$ = \biguplus _{m \le l_{1} < l_{2} \le n} X^{l_{1}}_{l_{2}} $

$\biguplus \{ ( p_{1} \cdots p_{m} \rightarrow \cdots \rightarrow c^{m} 
\rightarrow 
\sigma^{k_{1}}_{1} \rightarrow \cdots \rightarrow c^{k_{1}} \rightarrow 
 \sigma^{k_{2}}_{1} \rightarrow \cdots \rightarrow c^{k_{2}} \cdots 
 \rightarrow c^{n}  ) \ | \ |\sigma^{k_{i}}_{j}|=k_{i}    \}$.

Now for $m \le l_{1} < l_{2} \le n$, 

$Mob (X^{l_{1}}_{l_{2}}) = $

$\Sigma_{q_{1} \cdots q_{l_{1}} \neq c^{l_{1}}}
Mob (\{ (p_{1} \cdots p_{m} \rightarrow \theta_{1} \rightarrow \cdots
 \rightarrow \theta_{r} \rightarrow   q_{1} \cdots q_{l_{1}}) \ | \ 
 p_{1} \cdots p_{m} < \theta_{1} < \cdots  < 
 q_{1} \cdots q_{l_{1}}    \}) \cdot $

$Mob (\{ (q_{1} \cdots q_{l_{1}} \rightarrow \tau_{1} \rightarrow 
\cdots \rightarrow \tau_{s} \rightarrow c^{l_{2}}) \ | \ 
|\tau_{i}| = l_{2}      \}) \cdot $

$Mob ( \{ (c^{l_{2}} \rightarrow \sigma^{k_{1}}_{1} \rightarrow \cdots 
\rightarrow c^{k_{1}} \rightarrow \sigma^{k_{2}}_{1} \rightarrow \cdots 
\rightarrow c^{k_{2}} \rightarrow \cdots \rightarrow c^{n}) \ | \ 
|\sigma^{k_{i}}_{j}| = k_{i}        \})$.

By Lemma {\rmfamily \ref{import}} , we have 

$Mob (\{ (q_{1} \cdots q_{l_{1}} \rightarrow \tau_{1} \rightarrow 
\cdots \rightarrow \tau_{s} \rightarrow c^{l_{2}}) \ | \ 
|\tau_{i}| = l_{2}      \}) = 0$.

Hence we have 
$Mob (X^{l_{1}}_{l_{2}}) =0$ 

Therefore we have 
$\mu (p_{1}, \cdots p_{m}, c^{n}) = $

$
Mob (\{ ( p_{1}, \cdots p_{m} \rightarrow \cdots \rightarrow c^{m} 
\rightarrow 
\sigma^{k_{1}}_{1} \rightarrow \cdots \rightarrow c^{k_{1}} \rightarrow 
 \sigma^{k_{2}}_{1} \rightarrow \cdots \rightarrow c^{k_{2}} \cdots 
 \rightarrow c^{n}  ) \ | \ |\sigma^{k_{i}}_{j}|=k_{i}    \}) $

$= Mob (\{  (p_{1} \cdots p_{m} \rightarrow \theta_{1} 
 \cdots \theta_{r} \rightarrow c^{m})  \ | \ |\theta_{i}| = m    \}) \cdot $

$Mob (\{ (c^{m} 
\rightarrow 
\sigma^{k_{1}}_{1} \rightarrow \cdots \rightarrow c^{k_{1}} \rightarrow 
 \sigma^{k_{2}}_{1} \rightarrow \cdots \rightarrow c^{k_{2}} \cdots 
 \rightarrow c^{n}  ) \ | \ |\sigma^{k_{i}}_{j}|=k_{i}    \})$

$= (-1)^{m} \cdot \mu (c^{m},c^{n})$.

Hence we obtain the derived result.   
 \ \ \ \ \ \ \ \ \ \ $\Box$

\begin{lem}\label{example}

We have 

$\mu (\phi,c) = s-1  , \  \ \mu (a_{i},c) = -1,  \  
\mu (c,c^{2}) = 2s-1, \ \mu (a_{i},c^{2}) =-2s+1$.

\end{lem}

Now we put 
$T(k,n) := MobPat \{  c^{k} \rightarrow \theta_{1} \rightarrow \cdots 
\rightarrow \theta_{r} = c^{n} \ | \ |\theta_{i}| = n   \}$, $T(n,k) = 0$
 for $0 \le k < n$, and $T(n,n) :=-1$ for $0 \le n$ .

\begin{lem}

For $0 \le k \le n$, we have 

\[
T(k,n) = -\Sigma_{i=k}^{n}
(\begin{array}{c}
n \\
i
\end{array}
)
\cdot s^{n-i} \cdot (-1)^{n-i}.
\]

\end{lem}

\textbf{PROOF}
 
In case when $k = n$, it is trivial.

In case when $ 0 \le k < n$,  we have 

$T(k,n) = - \Sigma_{i=k}^{n} M((k,k),(n,i)) \cdot (-1)^{n-i}$.

Now we have 
\[
M((k,k),(n,i)) = 
(
\begin{array}{c}
n \\
i
\end{array}
)
\cdot s^{n-i}.
\]

Hence we obtain the derived result.

\begin{lem}

For $0 \le k \le l$,  we have 
$T(k,l) - T(k-1,l-1) = -sT(k,l-1)$. 

\end{lem}

\textbf{PROOF}

In case when $k=l$, it is trivial.

It is enough to show the case of $k<l$.
Then we have 

$T(k,l) - T(k-1,l-1)$

\[
= -\Sigma_{i=k}^{l}
(\begin{array}{c}
l \\ 
i
\end{array}
)
\cdot s^{l-i} \cdot (-1)^{l-i}  
+
\Sigma_{i=k-1}^{l-1}
(\begin{array}{c}
l-1 \\ 
i
\end{array}
)
\cdot s^{l-i-1} \cdot (-1)^{l-i-1}
\]

\[
=  -\Sigma_{i=k}^{l-1}
(\begin{array}{c}
l \\ 
i
\end{array}
)
\cdot s^{l-i} \cdot (-1)^{l-i} -1
+ 
\Sigma_{i=k-1}^{l-2}
(\begin{array}{c}
l-1 \\ 
i
\end{array}
)
\cdot s^{l-i-1} \cdot (-1)^{l-i-1} +1
\]

\[
= -\Sigma_{i=k}^{l-1}
\{ 
(\begin{array}{c}
l \\
i
\end{array})
-
(\begin{array}{c}
l-1 \\
i-1
\end{array})
\} \cdot s^{l-i} \cdot (-1)^{l-i}
\]

\[
=-s \cdot\Sigma_{i=k}^{l-1}
(\begin{array}{c}
l-1 \\
i
\end{array})
\cdot s^{l-i-1} \cdot(-1)^{l-i-1} = -sT(k,l-1).
\]

Hence we obtain the derived result.    \ \ \ \ \ \ \ \ \ \ $\Box$

\begin{lem}\label{c-eq}

Suppose $1 \le m < n$, then we have 

$\mu (c^{m}, c^{n}) = \Sigma_{k=m}^{n-1} \mu (c^{m},c^{k}) T(k,n)$.

\end{lem}

\textbf{PROOF}
 
 We heve 
$\mu (c^{m}, c^{n})$

$ = Mob (\{ (c^{m} \rightarrow \sigma^{m_{1}}_{1} \rightarrow \cdots 
\rightarrow c^{m_{1}}  \rightarrow \sigma^{m_{2}}_{1} \rightarrow 
\cdots \rightarrow c^{m_{2}} \rightarrow \cdots 
 c^{m_{r-1}} \rightarrow \sigma^{m_{r-1}}_{1} \cdots \rightarrow 
 c^{m_{r}} = c^{n}) \ | \ |\sigma^{m_{i}}_{j}| = m_{i} \})$.

In  case of $2 \le r$, we put $k = m_{r-1}$ and in case of 
$r=1$, we put $k =m$. By Lemma
{\rmfamily \ref{mobius}}, we obtain the derived result.

$\Box$

\begin{lem}\label{relation}
For $1 \le m <n$, we have 

$\mu (c^{m},c^{n}) - \mu (c^{m-1},c^{n-1}) = s \mu (c^{m},c^{n-1})$.

\end{lem}

\textbf{PROOF}

In  case of $n=2$, we see 
$\mu (c,c^{2}) - \mu (\phi,c) = (2s-1) -(s-1) =s , \ \mu (c,c)=1$ 
. Hence  this lemma holds for $n = 2$.
holds.

We prove by induction on $n$.
Suppose that the relation holds for $n-1$. Then we have 

$\mu (c^{m},c^{n}) - \mu (c^{m-1},c^{n-1}) $

$= \Sigma_{k=m}^{n-1} \mu (c^{m},c^{k}) T(k,n)
- \Sigma_{k=m-1}^{n-2} \mu (c^{m-1},c^{k}) T(k,n-1)$

$=\Sigma_{k=m}^{n-1} \mu (c^{m},c^{k}) T(k,n) - 
\mu (c^{m},c^{k}) T(k-1,n-1) + \mu (c^{m},c^{k}) T(k-1,n-1)$

$- \mu (c^{m-1},c^{k-1}) T(k-1,n-1)$

$= \Sigma_{k=m}^{n-1} \mu (c^{m},c^{k}) (-s) T(k,n-1)
+ s  \mu (c^{m},c^{k-1}) T(k-1,n-1)$

$=\Sigma_{k=m}^{n-1} \mu (c^{m},c^{k}) (-s) T(k,n-1) 
+s \Sigma_{k=m+1}^{n-1} \mu (c^{m},c^{k-1}) T(k-1,n-1)$

$=(-s) \mu (c^{m},c^{n-1}) T(n-1,n-1)$

$=s\mu (c^{m},c^{n-1}) $.

\begin{lem}\label{relation2}
For $1 \le m \le n$, we have 

$\mu (a_{1}^{m},c^{n}) + \mu (a_{1}^{m-1},c^{n-1}) = 
s \mu (a_{1}^{m},c^{n-1})$.

\end{lem}

\textbf{PROOF}

In case of $m < n$, we have by Lemma {\rmfamily \ref{relation}}.
In case of $m = n$, from $a_{1}^{m} \nleqq c^{n-1}$, 
$\mu (a_{1}^{m},c^{n}) = (-1)^{m}$ and 
$\mu (a_{1}^{m-1},c^{n-1})= (-1)^{m-1}$, therefore the right hand side $=0$.
 Hence  we obtan the derived result.

\begin{lem}

For $1 \le m \le n$,  $\mu (a_{1}^{m} , c^{n})$ 
is coefficient of $X^{n-m}$ in $T^{s}_{m+n}(X)$.

\end{lem}

\textbf{PROOF}

If 
$m+n = 2 , i.e  \ m=n=1$, we have 
$T^{s}_{2}(X) = s(s-1)X^{2} -1 , \ \mu (a_{1},c) = -1 $. 
Hence this lemma holds.

If 
$3 \le m+n$, by the relation 
$T^{s}_{k+2}(X) + T^{s}_{k}(X) = s X \cdot T^{s}_{k+1}(X)$ 
and Lemma {\rmfamily \ref{relation2}} we obtain the derived result.

\begin{lem}

 For $n \in \mathbb{N}$ we heve $\mu(\phi , c^{n}) = s^{n-1} (s-1)$.

\end{lem}

\textbf{PROOF}

If $n=1$, our claim follows from Lemma 
{\rmfamily \ref{example}}. We show by induction. We suppose that 
$\mu(\phi , c^{k}) = s^{k-1} (s-1)$ when $k \le n-1$

Now we have 

$\mu (\phi , c^{n})$

$= \Sigma_{k=1}^{n}  MobPat(\phi , c^{k}) \mu(c^{k}, c^{n})$
\ \ \ \ \ \ \ \ \ \ \ \ \ (by Lemma{\rmfamily \ref{import}})

$=\Sigma_{k=1}^{n} \bigl( -\Sigma_{i=0}^{k} s^{i} 
M((0,0),(k,k-i)) (-1)^{i}  \bigr) \mu (c^{k},c^{n})$

\[
=\Sigma_{k=1}^{n} \bigl( -\Sigma_{i=0}^{k} s^{i}
(
\begin{array}{c}
k \\
i 
\end{array}
)
(-1)^{i}  \bigr) \mu (c^{k},c^{n})
\ \ \ \ \ \ \ \ \ \ \ \ \ \ \ \ \ \ \ \ \ \ \ \ \ \ \ \ \ \ \ \ \ \ \ \ \ \ \ 
\ \ \ \ \ \ \ \ \ \ \ \ \ \ \ \ \ \ \ \ \ \ \ \ \ \ \ \ \ \ \ \ \ \ \ \ \ \ \ 
\ \ \ \ \ \ \ \ \ \ \ \ \ \ \ \ \ \ \ \ \ \ \ \ \ \ \ \ \ \ \
\]

$=\Sigma_{k=1}^{n} -(1-s)^{k} \mu (c^{k},c^{n})$.

And we have 

$\mu (\phi,c^{n}) - s \mu (\phi, c^{n-1})$

$=\Sigma_{k=1}^{n} -(1-s)^{k} \mu (c^{k},c^{n})
-\Sigma_{k=1}^{n-1} -(1-s)^{k} \mu (c^{k},c^{n-1})$

$= \Sigma_{k=1}^{n-1} -(1-s)^{k} \bigl(  \mu (c^{k},c^{n})  
-s \mu (c^{k},c^{n-1})    \bigr)  - (1-s)^{n}$

$=\Sigma_{k=1}^{n-1} -(1-s)^{k} \mu (c^{k-1},c^{n-1}) -(1-s)^{n}$
\ \ \ \ \ \ \ \ \ \ \ (by Lemma{\rmfamily \ref{c-eq}})

$=-(1-s) \bigl(    \Sigma_{k=1}^{n-1} (1-s)^{k-1} \mu (c^{k-1},c^{n-1})   
\bigr) -(1-s)^{n}$

$=-(1-s) \bigl(    \Sigma_{k=0}^{n-2} (1-s)^{k} \mu (c^{k},c^{n-1})   
\bigr) -(1-s)^{n}$

$=-(1-s)  
\bigl(  -\mu (\phi, c^{n-1}) - \mu(\phi, c^{n-1}) - (1-s)^{n-1})      \bigr)
-(1-s)^{n}$

$=0$.

So we have $\mu (\phi,c^{n}) = s \mu (\phi, c^{n-1})$. Hence we obtain 
the derived result.        $\Box$

\begin{lem}

For $n \in \mathbb{N}$,   $\mu (\phi,c^{n})$ 
is the coefficient of $X^{n}$ in $T^{s}_{n}(X)$.

\end{lem}

Therefore we have the following theorem.

\begin{thm}\label{last}

For $0 \le m \le n$ , 
$\mu (a_{1}^{m},c^{n})$ is the coefficient of $X^{n-m}$ in 
$T^{s}_{m+n}(x)$. 

\end{thm}

\begin{cor}

Conjecture {\rmfamily \ref{sagan-conj}} is true.

\end{cor}

{\large \textbf{ACKNOWLEDGEMENT}}

The author wishes to thank Professor Jun Morita for his valuable advice. And 
he is  also grateful to Professor Daisuke Sagaki, Sho Matsumoto for their
 helpful comments.

\renewcommand{\refname}{REFERENCE}

\end{document}